\documentclass[12pt]{article}

\usepackage{amssymb}

\newtheorem{theorem}{Theorem}[section]
\newtheorem{lemma}[theorem]{Lemma}
\newtheorem{corollary}[theorem]{Corollary}
\newtheorem{proposition}[theorem]{Proposition}

\newtheorem{definition}{Definition}
\newtheorem{claim}{Claim}

\newcommand\qed{\hfill\hbox{\vrule height2mm width2mm}\medskip\par}

\newenvironment{proof}{\begin{trivlist}  \item[\hskip\labelsep{\it Proof.}]}%
                     {\hfill {\qed}\end{trivlist}}

\newcommand\sub{\subseteq}

\newcommand{\N}{\mathbb{N}}
\newcommand{\R}{\mathbb{R}}
\newcommand{\Z}{\mathbb{Z}}

\newcommand\dts{,\dots ,}

\newcommand\be{\begin{equation}}
\newcommand\bel[1]{\begin{equation}\label{#1}}
\newcommand\ee{\end{equation}}
\newcommand\ba[1]{\begin{array}{#1}}
\newcommand\ea{\end{array}}
\newcommand\bl{\begin{lemma}}
\newcommand\bll[1]{\begin{lemma}\label{#1}}
\newcommand\el{\end{lemma}}
\newcommand\bi{\begin{itemize}}
\newcommand\ei{\end{itemize}}
\newcommand\bt{\begin{theorem}}
\newcommand\btl[1]{\begin{theorem}\label{#1}}
\newcommand\et{\end{theorem}}
\newcommand\ben{\begin{enumerate}}
\newcommand\een{\end{enumerate}}
\newcommand\bpr{\begin{proposition}}
\newcommand\bprl[1]{\begin{proposition}\label{#1}}
\newcommand\epr{\end{proposition}}
\newcommand\bcl{\begin{claim}}
\newcommand\ecl{\end{claim}}
\newcommand\bprf{\begin{proof}}
\newcommand\eprf{\end{proof}}
\newcommand\bdf{\begin{definition}}
\newcommand\edf{\end{definition}}

\title{Randomness, pseudorandomness and models of arithmetic\\
~~\\
{\it\Large  dedicated to Alan Woods}}
\author{Pavel Pudl\'ak%
\thanks{Partially supported by grants IAA100190902 of GA AV \v{C}R and
  GAP202/12/G061 of GA \v{C}R.}}

\begin{document}

\maketitle

\begin{abstract}
  Pseudorandmness plays an important role in number theory, complexity
  theory and cryptography. Our aim is to use models of arithmetic to
  explain pseudorandomness by randomness. To this end we construct a
  set of models $\cal M$, a common element $\iota$ of these models and
  a probability distribution on $\cal M$, such that for every
  pseudorandom sequence $s$, the probability that $s(\iota)=1$ holds
  true in a random model from $\cal M$ is equal to $1/2$.
\end{abstract}

\section{Introduction}

A pseudorandom sequence is an infinite sequence of $-1$s and $1$s
computable in nondeterministic polynomial time that is not correlated
with any polynomial time computable function. Such sequences can also
be viewed as sets in $\bf NP\cap coNP$; thus we can also talk about
pseudorandom sets. Intuitively, a pseudorandom set splits every set in
$\bf P$ into two sets of equal density. There are some natural and
important candidates for pseudorandom sequences in number theory such
as the Liuville function (closely related to the M\"obius function).

Our main result is a construction of a set of models $\cal M$, a
common element $\iota$ of these models and a probability distribution
on $\cal M$, such that for every pseudorandom sequence $s$, the
probability that $s(\iota)=1$ holds true in a random model from $\cal
M$ is equal to $1/2$. Thus pseudorandomness of sequences manifests
itself in $\cal M$ as genuine randomness. Admittedly, this result is
weak, because it concerns only one common element of the models. We
present it as a proof of the concept that results of this kind are
possible. We suggest some ways of extending this result in
Section~\ref{remarks}.

We prove our result by using restricted ultrapowers, which are
ultrapowers in which the sets of the ultrafilter and the  functions are
elements suitable classes. The history of restricted ultrapowers goes
back to Skolem (see~\cite{kochen-kripke}). We will start with a model
$M_0$, a core of our construction, constructed from the set of all
polynomial time computable functions reduced by a suitable ultrafilter
on the complexity class $\bf P$. By extending the
class of functions and the ultrafilter in various ways, we obtain a
set of models $\cal M$ extending the core model. Remarkably, there is a natural
way of defining a probability measure on $\cal M$.

In the last section we present some philosophical speculations about
the nature of pseudorandomness.

\section{Preliminaries}

\subsection{Random sequences}

We will study sequences $s:\N\to\{\pm 1\}$. Let $\cal S$ denote the
set of all such sequences with the uniform distribution. 
Intuitively, $Pr[s(n)=1]=Pr[s(n)=-1]=\frac 12$ and these events are
independent for different numbers.
Formally, there is a Lebesgue measure $m$ on $\cal S$ that is
uniquely determined by 
\[
m(\{s\in{\cal S};\ s(1)=a_1\wedge\dots\wedge s(n)=a_n \})=2^{-n}
\]
for all $n$ and all strings $a$ of $\pm 1$s.  We say that a {\em
  ``random sequence satisfies $P$,''} if the probability that a random
sequence satisfies $P$ is 1. If a random sequence satisfies
$P_1,P_2,\dots$, then it also satisfies $\bigwedge_i P_i$.

The basic fact about random sequences is the Law of Large Numbers
\[
Pr\left[\lim_{n\to\infty}{\sum_{i=0}^{n-1}s(i)}/n=0\right]=1.
\]
This theorem, however, does not provide information about the rate of
convergence. Much more precise theorems have been proven, in
particular Khinchin's Law of Iterated Logarithm
\[
Pr\left[\limsup_{n \to \infty} {\sum_{i=0}^{n-1}s(i)}/{\sqrt{2n \ln\ln n}} = 1\right]=1, 
\]
which, in particular, implies that for every $\epsilon>0$,
\bel{eRH}
Pr\left[\lim_{n\to\infty}{\sum_{i=0}^{n-1}s(i)}/{n^{\frac 12+\epsilon}}=0\right]=1.
\ee

\subsection{Algorithmic randomness}

We want to formalize the concept that a sequence $s$ ``looks like a
random sequence''. Here the sequence $s$ is fixed, so we cannot use
probability theory. The basic idea is that $s$ must satisfy many
properties that random sequences satisfy with probability 1. E.g., we
certainly want the property used in the Law of Large Numbers. Further,
we want to consider properties that can be algorithmically
decided. Therefore it is more natural to talk about satisfying {\em tests}
instead of properties. 

The study of such concepts has long history and many researchers
contributed to it, including Kolmogorov, Chaitin, Levin, Schnorr and
Martin-L\"of. This research area is called {\em algorithmic
  randomness}. We mention one of the concepts that is studied there, so
that we can compare it with the concept that we will introduce.

A {\em martingale} is a function $F:\{\pm 1\}^*\to\R^+$ such that
\[
F(a_1\dts a_n)=\frac 12(F(a_1\dts a_n,-1)+F(a_1\dts a_n,+1))
\]

\begin{definition}[Schnorr \cite{schnorr}]
  A sequence $s:\N\to\{\pm 1\}$ is {\em {\bf P}-random,} if for every
  polynomial time computable martingale
\[
\limsup_{n\to\infty}F(s(0),s(1)\dts s(n-1))<\infty.
\]
\end{definition}

Another concept relevant to this paper is the concept of {\em
  pseudorandom number generator}, which will be abbreviated by
{PRG}. A PRG is an algorithm to produce a long string of numbers, usually
just 0s and 1s, from a short random string, called the {\em seed}. So
in this case we do not have only one infinite sequence, but a small set of
finite strings with a probability distribution. This concept plays an
important role in the theories that study computational complexity and
cryptography. 

There is one important difference in how the computational resources
are bounded in the two mentioned approaches. When testing a sequence
$s$ for {\bf P}-randomness, the testing algorithm receives the whole
string $(s(0),s(1)\dts s(n-1))$ as an input and can use time
polynomial in $n$. In contrast, algorithms testing pseudorandomness
get $n$ represented in binary and can use time polynomial in the
length of $n$, i.e., they run in time polynomial in $\log n$.

\subsection{Pseudorandom sequences}

Our concept is closer to the theory of pseudorandomness than to
algorithmic randomness. That is why we use the word {\em
  pseudorandom}. The reader should, however, be cautioned that there
are concepts with similar names that differ significantly from ours.

\begin{definition}\label{def2}
A sequence $s:\N\to\{\pm 1\}$ will be called {\em pseudorandom} if
\ben
\item $s$ is computable in nondeterministic polynomial time, 
\item for every polynomial time computable function  $f:\N\to\{\pm 1\}$,
\bel{e0}
\lim_{n\to\infty}\sum_{i=0}^{n-1} f(i)s(i)/n =0.
\ee
\een
\end{definition}
Condition 1. means that although we may be not able to compute the
value $s(n)$ in polynomial time, if somebody gives us a ``witness'' we
are able to check the correct value of $s(n)$ in polynomial
time. (There should always be witnesses for either the value $1$ or
the value $-1$, but never for both.) We can identify $s$ with the set
$\{n;\ s(n)=1\}\in\bf NP\cap coNP$, so we can also talk about {\em
  pseudorandom sets.}

Condition 2. means that $s$ is little correlated with any sequence
computable in polynomial time. One can consider various modifications
of this condition. For instance, one can allow stronger tests, say, functions
computable in appropriately defined subexponential time. One can also
impose stronger bounds on correlation. If, for example, we required
that 
\[
\lim_{n\to\infty}\sum_{i=0}^{n-1} f(i)s(i)/n^\alpha =0,
\]
for some $1/2<\alpha<1$, then the correlation would be exponentially small
(on finite initial segments; recall that the input size is $\log n$).

\medskip Another motivation for the definition above is the {\em
  M\"obius Randomness Principle} proposed by Peter Sarnak
(see~\cite{green}). According to this principle the M\"obius function
$\mu$ is not correlated to any ``low-complexity'' function
$F:\N\to[-1,1]$ in the sense that
\[
\lim_{n\to\infty}\sum_{i=1}^{n} F(i)\mu(i)/n =0.
\]
This is like saying that $\mu$ is pseudorandom, except for a few 
minor differences. First, $\mu$ takes on not only the values $\pm 1$, but
also $0$. This does not seem very important for studying the concept
of pseudorandomness (see Proposition~\ref{mu-lambda} below).  Second,
it is not specified what ``low complexity'' means. This leaves open
the possibility of studying various specific versions of the conjecture.
Third, the tests
are functions $F$ whose range is in the whole interval $[-1,1]$. We
will show below that this is also an irrelevant difference.

\bpr
Suppose that $s:\N\to\{\pm 1\}$ is pseudorandom. Let $F:\N\to[-1,1]$
be a polynomial time computable function whose values are binary
rationals. Then
\[
\lim_{n\to\infty}\sum_{i=1}^{n} F(i)s(i)/n =0.
\]
\epr
\bprf
Let $s$ and $F$ be given.
One can easily show that $s$ is pseudorandom also with respect to
polynomial time computable functions $f:\N\to\{-1,0,1\}$. (Hint: write 
$f=\frac 12 f_+ + \frac 12 f_-$, where $f_+(n)=1$ if $f(n)=1$,
otherwise $f_+(n)=-1$, and $f_-(n)=-1$ if $f(n)=-1$,
otherwise $f_-(n)=1$.)

Represent $F$
as a weighted sum of such functions, 
$
F(n)=\sum_{j=0}^\infty 2^{-j}f_j(n).
$
Then
\[
\lim_{n\to\infty}\sum_{i=1}^{n} F(i)s(i)/n 
=\lim_{n\to\infty}\sum_{i=1}^{n} \sum_{j=0}^\infty 2^{-j}f_j(n)s(i)/n
=\sum_{j=0}^\infty 2^{-j}\lim_{n\to\infty}\sum_{i=1}^{n} f_j(n)s(i)/n
=0
\]
because the infinite sum converges absolutely.
\eprf

Say that a set of numbers $A$ has {\em positive density,} if
$\liminf_{n\to\infty}{|A\cap[0,n-1]|}/n>0$. 
\begin{corollary}\label{cor1}
If $X$ is a pseudorandom set, then neither $X$ nor its complement
contains a set $A\in\bf P$ of positive density.
\end{corollary}

A possible way of stating the M\"obius Randomness Principle is to say
that $\mu$ is pseudorandom in the sense of our definition (extended to
sequences of $-1$, $0$ and $1$). 
A closely related function is the Liouville function $\lambda$. It is
defined by $\lambda(n)=(-1)^k$, where $k$ is the number of prime
factors of $n$ counted with their multiplicity. 
\bpr\label{mu-lambda}
$\mu$ is pseudorandom if and only if $\lambda$ is.
\epr
\begin{proof}
1. Suppose $\mu$ is pseudorandom. Let a polynomial time computable
function  $f:\N\to\{\pm 1\}$ and $\epsilon>0$ be given. Let $n_0$ be
such that $\sum_{k>n_0}k^{-2}<\epsilon$. We have
\[
|\lim_{{n\to\infty}}\sum_{i=1}^{n} f(i)\lambda(i)/n|\leq
\]\[
\sum_{1\leq k\leq n_0}|\lim_{n\to\infty}
\sum_{\stackrel{1\leq k^2i\leq n,} {i \mbox{ \tiny square free}}} 
f(k^2i)\lambda(k^2i)/n\ |+
\sum_{k>n_0}|\lim_{n\to\infty}\sum_{\stackrel{1\leq k^2i\leq n,}{ i\mbox{ \tiny
      square free}}} f(k^2i)\lambda(k^2i)/n\ | <
\]\[
\sum_{1\leq k\leq n_0}|\lim_{n\to\infty}\sum_{1\leq k^2i\leq n}
  f(k^2i)\mu(i)/n\ |+ \epsilon = \epsilon.
\]

2. Now suppose that $\mu$ is not pseudorandom.  Let  $f:\N\to\{\pm
1\}$ be a polynomial time computable
function  and $\epsilon>0$ such that 
\mbox{$\lim_{n\to\infty}\sum_{i=1}^{n} f(i)\mu(i)/n=\epsilon.$}
Let $n_0$ be as above. 
Furthermore, we can suppose that 
\mbox{$\lim_{n\to\infty}\sum_{i=1}^{n} f(i)\lambda(i)/n=0,$}
because
otherwise we would be done. In a similar fashion as above,
decompose the $\sum_{i=1}^{n} f(i)\lambda(i)/n$ into three terms
\ben
\item the sum over square free numbers $i$,
\item the sum over numbers $i$ divisible by $k^2$ for some $k\leq n_0$, and
\item the sum over the remaining numbers $i$.
\een
By our assumptions, the limit of the first sum is $\epsilon$, the
sum of the limits of the sums 2. and 3. is $-\epsilon$, the sum 3. is
$>-\epsilon$. Hence, if we define 
\[
g(n)=\left\{
\begin{array}{rl}
-f(n)&\mbox{ if $n$ is divisible by $k^2$ for some }k\leq n_0,\\
f(n)&\mbox{ otherwise,}
\end{array}\right.
\]
we obtain
$\lim_{n\to\infty}\sum_{i=1}^{n} g(i)\lambda(i)/n>0.$
\end{proof}

Some special cases of the M\"obius randomness principle have been
proven. The first one was the Prime Number Theorem, which is the case
of $f\equiv 1$. (The question
whether the bound on correlation can be improved to the form (\ref{eRH}) is
the Riemann Hypothesis.)
Recently B. Green
proved the principle for $AC^0$, \cite{green}. Let us say that a sequence is
$AC^0$-pseudorandom if it satisfies Definition~\ref{def2} with the
condition 2. weakened to $AC^0$-computable. Then one can state Green's
result as follows.
\bt[\cite{green}]
The M\"obius function is $AC^0$-pseudorandom.
\et
The Liouville function is also $AC^0$-pseudorandom.

It will be very difficult to prove that some sequence is pseudorandom,
because the existence of pseudorandom sequences implies $\bf P\neq
NP$. For specific functions, it may be even harder.  If the M\"obius
function is pseudorandom, then integers cannot be factored in
polynomial time.

In the opposite direction, we know that hardness of factoring implies
the existence of pseudorandom sequences (we are not able to prove that
it implies the pseudorandomness of the M\"obius function,
though). This is because 
\ben
\item there are constructions of permutations that are
{\em one-way functions} provided that factoring is hard, 
\item there is a construction of a {\em hard-core predicate} from any
  one-way permutation, and
\item hard-core predicates are very closely related to pseudorandom sequences. 
\een

We will now explain this connection in more detail, but for the sake
of brevity, we will skip the definition of a one-way function. Let
$1\leq k_1<k_2<\dots$ be a sequence of integers that grow at most
polynomially, and let $F_j:\{0,1\}^{k_j}\to\{0,1\}^{k_j}$, $j=1,2,\dots,$
be a sequence of permutations. Suppose that these numbers and
functions are uniformly computable in polynomial time. We say that
functions $B_j:\{0,1\}^{k_j}\to\{0,1\}$, $j=1,2,\dots,$ are hard-core
predicates for the functions $F_j$, if $B_j$ are uniformly computable
in polynomial time, and for every function $g(x)$ computable by a
randomized polynomial time algorithm
\bel{e05}
Pr[g(F_j(x))=B_j(x)]=\frac 12 \pm \frac 1{k_j^{\omega(1)}},
\ee
where the probability is taken over uniformly distributed
$x\in\{0,1\}^{k_j}$ and random bits of the algorithm for $g$; further,  
$\omega(1)$ is the standard notation for functions going to
infinity. In plain words, this means that $B_j(x)$ can be
predicted from  $F_j(x)$ only with negligible probability (which, in
particular, 
implies that it is difficult to invert $F_j$). 

The concept that is closely related to pseudorandom sequences (as
defined in this paper) is the sequence $B_j(F^{-1}(y))$, $j=1,2,\dots$. To get
a pseudorandom sequence $s$, we only need to connect the bits
$B_j(F^{-1}(y))$ into one infinite sequence:
\[
s(n)=(-1)^{B_j(F^{-1}(n-\sum_{i<j}2^{k_i}))},
\]
where $\sum_{i<j}2^{k_i}\leq n<\sum_{i\leq j}2^{k_i}$ and 
where we are identifying $\{0,1\}^{k_j}$ with $0,1,\dts 2^{k_j-1}-1$.

In order to show that $s(n)$ is pseudorandom, we have to address only
one small complication. While in (\ref{e0}) of the definition of
pseudorandomness we consider all initial segments, in (\ref{e05}) of
the definition of the hard-core predicate we only consider correlation
over the entire interval $[0,2^{k_j}-1]$ (but we have better
convergence). We need to show that the correlation of $B_j(F^{-1}(y))$
with polynomial time functions is also low on an initial segments
$[0,a]$ of $[0,2^{k_j}-1]$.

Suppose $g$ has positive correlation with $B_j(F^{-1}(y))$ on
$[0,a]$.  If we knew $a$, we could define $g'$ that has positive
correlation with 
$B_j(F^{-1}(y))$ on the entire interval
by putting $g'(x)=g(x)$ on $[0,a]$ and $g'(x)=1-g(x)$ on the
rest. Since we cannot assume that we know $a$, we have to do something
slightly more complicated. Note that we are actually assuming that
there are infinitely many indices $j$ for which $g$ has positive
correlation with 
$B_j(F^{-1}(y))$ on some interval
$[0,a_j]$. The ratios $a_j/2^{k_j}$ have some limit
point $\alpha$, $0<\alpha<1$. Take a rational number $\beta$ close to
$\alpha$ (or $\alpha$ itself if it is rational). Then use $\lceil\beta
2^{k_j}\rceil$ as switching points.

This finishes a sketch of the proof of the following proposition.
\bpr 
If there exists no probabilistic polynomial time algorithm for
factoring integers, then there exists a pseudorandom sequence.  
\epr

{\it Remarks.} 1. The construction actually gives a concrete sequence,
but since its definition is rather complicated, we do not present it
here.

2. We have not fully used the assumption about factoring; one can show
pseudorandomness of the constructed sequence in a little stronger sense.

3. For the concepts and results used above, see~\cite{goldreich}.

\subsection{Theories}

We need a theory in which it is possible to formalize polynomial time
computations. A natural theory in which this is possible is Cook's
$PV$, \cite{cook,krajicek}. This theory has function symbols for all
polynomial time computable function. The function symbols correspond
to algorithms based on recursion on
notation. Our result is quite general and, as such, holds for the stronger 
theory $PV^{\N}$ defined below.
\begin{definition}
~

\ben
\item $PV^{\N}$ is the theory axiomatized by all true \emph{universal sentences}
  in the language of $PV$.
\item $AR^{\N}$ is the theory consisting of \emph{all true sentences} in the
  language of $PV$.
\een
\end{definition}
The theory $PV^{\N}$ is a conservative extension of the theory of all
true $\Pi^0_1$ arithmetical sentences plus the axiom $\forall x\exists
y\ y=x^{\lceil\log(x+1)\rceil}$ (this axiom guarantees that the provably
total functions grow sufficiently fast and thus enable us to define
polynomial time computations). The theory $AR^{\N}$ is essentially 
{\em True Arithmetic,} except that we use the richer language of $PV$.


We focus on the complexity class ${\bf P}$, as this is the most
interesting case, but in fact the same result can be proven for
concepts based on other classes. An interesting case is the class $\bf
AC^0$ because of the result of Green mentioned above. The
theory corresponding to this class is $V^0$, see \cite{cook-nguyen}.

\subsection{Random models}

Our aim is to represent pseudorandomness in a different way. The basic
idea is to 
study this concept using a set of nonstandard models equipped with a
probability distribution. Let $\sigma$ be a first order formula in the
language of $PV$ defining a sequence $s\in\cal S$, and suppose that we
have a set of models $\cal M$ with a probability distribution $\nu$.
Let $\phi$ be a first order formula. Then we can say {\em `$s$
  satisfies the property $\phi$ with probability $p$'} if the
probability that in a random model from $\cal M$ the sequence defined
by $\sigma$ satisfies $\phi$ is $p$.

If we want to use formulas $\phi$ with parameters, i.e., free
variables for elements of models, we need to impose some structure on
$\cal M$. In this paper we will only consider the following
structure. There is one distinguished model $M_0$ such all other
models are its extensions. This enables us to speak about properties
parameterized by elements of $M_0$.

In this paper we say that a model $N$ is an
{\em extension} of a model $M$, if $M$ is a substructure of $N$.

In general the structure defined on $\cal M$ can be more
complicated. We can use various frames, like in the Kripke
semantics. If $N$ is one of the ``alternative worlds'' of $M$, then
$N$ should be an extension of $M$ (not necessarily proper).

An alternative approach is to use a Boolean valued model $M$ with a
boolean algebra $\cal B$ equipped with a probability measure, an
approach studied in~\cite{krajicek11}. This is, however, not
fundamentally different from the approach sketched above. Having such
a model, we can construct a set of models by taking all ultrafilters
on $\cal B$. Vice versa, having $\cal M$ and $\nu$ as above, we can
take the Boolean algebra of measurable subsets of $\cal M$ and define
a measure on this algebra in a natural way.

\section{The result}

\bt\label{t-main} 
There exists a model $M_0$ of $PV^{\N}$, an element $\iota\in M_0$,
a set $\cal M$ of models of $AR^{\N}$ and a probability measure
$\nu$ on a sigma algebra $\cal B$ of subsets of $\cal M$ such that 
\ben
\item models of $\cal M$ are extensions of $M_0$,
\item for every $PV$ formula $\phi(x_1\dts x_k)$ and every string of
  elements $a_1\dts a_k\in M_0$, the set $\{ M\in{\cal M};\
  M\models\phi(a_1\dts a_k)\}$ is in $\cal B$,
\item for every (definition of a) pseudorandom sequence $s$,
\[
Pr_\nu[M\models s(\iota)=1]=\frac 12.
\]
\een
\et

\begin{proof}
Let
  $K\sub\N$ be an infinite set. We define {\em density of sets of
    numbers with respect to $K$,} or {\em $K$-density,} to be the
  partial function defined by
\bel{e1}
dens_KX=\lim_{k\in K, k\to\infty}\frac{|X\cap[0,k-1]|}k.
\ee
($dens_KX$ is undefined when the limit does not exist.)

The proofs of the following easy facts are left to the reader.
\ben
\item $dens_K$ is finitely additive.
\item If $dens_KX=dens_KY=1$, then $dens_KX\cap Y=1$.
\item If $C$ is a countable set of sets of numbers, then there exists
  an infinite $K$ such that $dens_KX$ is defined (i.e., the limit
  (\ref{e1}) exists) for every $X\in C$.  
\een 
Let $\bf AR$ denote arithmetically definable sets of natural
numbers (which is the same as sets first-order definable in the
language of $PV$).  Let $K$ be an infinite set of numbers such that
$dens_KX$ is defined for every $X\in\bf AR$. 

The following fact is also easy.
\ben
\item[4.] Let $X\in\bf P$ and $Z$ be a pseudorandom set. Then $dens_K X\cap
  Z=\frac 12 dens_KX$.
\een

Let ${\cal F}_0$ be the filter in $\bf P$ consisting of all sets of
$K$-density 1. Let ${\cal U}_0$ be an ultrafilter extending ${\cal
  F}_0$. Hence all sets in ${\cal U}_0$ have positive density. 
Let $\bf FP$ denote the set of functions computable in polynomial time.
We define
$M_0$ to be the ultrapower constructed by taking $\bf FP$ modulo ${\cal U}_0$,
\[
M_0={\bf FP}/{\cal U}_0,
\]
with the $PV$ function symbols interpreted in the natural way. The
fact in $M_0$ all true universal $PV$ sentences are satisfied is an
immediate consequence of {\L}o\'s's theorem. 
The distinguished element $\iota\in M_0$ is defined to be the element of
${\bf FP}/{\cal U}_0$ represented by the identity function $id$ on
$\N$, in symbols $\iota=[id]_{{\cal U}_0}$.

\medskip Let $u_0=\{ U_1\supset U_2\supset\dots\}$ be a cofinal chain
in ${\cal U}_0$. We define the {\em density of a set with respect to
  $u_0$,} or {\em $u_0$-density,} to be the partial function defined
by%
\footnote{More precisely, we should also use the index $K$, but
  there is no danger of confusion, since $K$ is fixed for the rest of
  the proof.} 
\bel{e2} 
dens_{u_0}X =  \lim_{n\to\infty}\frac{dens_{K} X\cap U_n}{dens_K U_n}.  
\ee 
Again, one can easily prove that we can pick $u_0$ so that 
$u_0$-density is defined for every $X\in\bf AR$.
The following facts are immediate corollaries of 1. and 4.:
\ben
\item[5.] $dens_{u_0}$ is finitely additive.
\item[6.] $dens_{u_0}Z=\frac 12$ for every pseudorandom set $Z$.
\een
The set of models that extend $M_0$ will be defined using the
following set of ultrafilters on the boolean
  algebra $\bf AR$.
\[
\Omega = \{{\cal U};\ {\cal U}\mbox{ ultrafilter on }{\bf AR},\
{\cal U}_0\sub{\cal U}
\mbox{ and }\forall V\in{\cal U}\ dens_{u_0}V>0\}.
\]
\ben
\item[7.] If ${\cal F}\sub\bf AR$ is a filter such that
  $dens_{u_0}U>0$ for all $U\in\cal F$, then $\cal F$ can be extended
  to an ultrafilter belonging to $\Omega$.
\een
Let $\bf FAR$ denote arithmetically definable functions. Define
\[
{\cal M} = \{ M;\ M={\bf FAR}/{\cal U},\ {\cal U}\in\Omega\}.
\]
Note that for ${\cal U}_1\neq{\cal U}_2$, the ultrapower models are
different (they may be isomorphic, though). 
Since $\bf FP\sub FAR$ and ${\cal U}_0\sub{\cal U}$, for every
${\cal U}\in\Omega$, we have:
\ben
\item[8.] Every $M\in\cal M$ is an extension of $M_0$. 
\een
The fact that these models are models of True Arithmetic, is a
well-known consequence of {\L}o\'s's theorem.

For $X\in\bf AR$, let
\[
\Omega[X]=\{{\cal U};\ X\in\cal U\},
\]
and put 
\[
{\cal A}_0=\{\Omega[X];\ X\in\bf AR\}
\]
\ben
\item[9.] ${\cal A}_0$ is a Boolean algebra.
\een
\bl
$\Omega[X]=\Omega[Y]$ if and only if $dens_{u_0}X\triangle Y=0$ (where
$\triangle$ denotes the symmetric difference).
\el
\begin{proof} $dens_{u_0}X\setminus Y>0$. 
Let $\cal F$ be the filter in $\bf AR$ generated by $X\setminus Y>0$. 
By 7., $\cal F$ has an extension to ${\cal U}\in\Omega$. Hence  ${\cal
  U}\in\Omega[X]\setminus\Omega[Y]$. 
This gives us the forward implication.

Now suppose ${\cal U}\in\Omega[X]\setminus\Omega[Y]$, for some $\cal
U$. Then $X\setminus Y\in\cal U$. Since ultrafilters in $\Omega$ do
not contain sets of $u_0$-density 0, we have  $dens_{u_0}X\setminus
Y>0$. 
\end{proof}


This lemma enables us to define an additive measure $\nu_0$ on ${\cal
  A}_0$ by putting 
\[
\nu_0(\Omega[X])=dens_{u_0}X.
\]
In particular, $\nu_0(\Omega)=1$. 
\bl
If $A_1,A_2,\dots\in{\cal A}_0$ are pairwise disjoint 
and $\bigcup_n A_n\in{\cal A}_0$, then 
$\nu_0(A_n)=0$ for all $n$ except for a finite number of them.  
\el
\begin{proof}
To prove the claim, suppose the contrary. Let $X,X_1,X_2,\dots \in \bf
AR$ be such that $A_n=\Omega[X_n]$, for $n=1,2,\dots$, and
$\bigcup_n A_n=\Omega[X]$. We observe that
$dens_{u_0}X_i\cap X_j=0$ for $i\neq j$, because $A_i\cap
A_j=\emptyset$.

Let $Y_n=X\setminus\bigcup_{i=1}^{n-1}X_i$. We will show
that $dense_{u_0}Y_n>0$ for all $n$. Suppose that for some $n$,
$dense_{u_0}Y_n=0$. Let $m\geq n$ such that $dens_{u_0}X_m>0$. Since
$dens_{u_0}X_m\cap \bigcup_{i=1}^{n-1}X_i=0$, we have
$dens_{u_0}X_m\cap X=0$. This implies that
$\Omega[X_m]\cap\Omega[X]=\emptyset$. But this is impossible, because
$\Omega[X_m]\neq\emptyset$. Thus $dense_{u_0}Y_n>0$ for all $n$.

Extend the filter $\{ Y_n;\ n=1,2,\dots\}$ to an ultrafilter ${\cal
  U}\in\Omega$. Clearly ${\cal U}\in\Omega[X]$, but for no $n$, ${\cal
  U}\in\Omega[X_n]$. 
\end{proof}

An immediate corollary is:
\ben
\item[10.] $\nu_0$ is $\sigma$-additive.  
\een 
According to a basic theorem about extensions of measures, we can
extend the $\sigma$-additive probability measure $\nu_0$ defined on
the Boolean algebra ${\cal A}_0$ to a $\sigma$-additive probability
measure $\nu_1$ defined on a $\sigma$-algebra ${\cal A}_1$. Using the
bijection ${\cal U}\mapsto{\bf FAR}/{\cal U}$ we translate the measure
$\nu_1$ defined on a $\sigma$-algebra ${\cal A}_1$ to a measure $\nu$
defined on a $\sigma$-algebra $\cal B$ of subsets of $\cal M$.

\bigskip 
We will now prove the second condition of the theorem. 
Let $\phi(x_1\dts x_k)$ be a $PV$ formula, $a_1\dts a_k\in
M_0$, let ${\cal U}\in\Omega$ and let $M={\bf FAR}/{\cal
  U}$. Further, let $f_1\dts f_k\in\bf FAR$ be the functions
representing $a_1\dts a_k$ (in symbols, 
$a_i=[f_i]_{{\cal U}_0}$).

According to {\L}o\'{s}'s theorem, $M\models \phi(a_1\dts a_k)$
if and only if 
\[
\{ n;\ \N\models\phi(f_1(n)\dts f_k(n))\}\in \cal U.
\]
Hence the set of ultrafilters for which the models satisfy
$\phi(a_1\dts a_k)$ has the form $\Omega[X]$, for $X\in\bf
AR$. Therefore 
\[
\{ M;\ M\models \phi(a_1\dts a_k)\}\in\cal B.
\]

\medskip
It remains to prove the third condition of the theorem. Let $s$ be a
pseudorandom sequence and let $\psi(x)$ be a formula defining
$s(x)=1$. By {\L}o\'{s}'s theorem, the set of models satisfying
$\psi(\iota)$ corresponds to the set of ultrafilters such that 
$\{ n;\ \N\models\psi(n)\}\in \cal U$ (recall that $\iota=[id]_{{\cal
    U}_0}$). Thus we have
\[
\begin{array}{rcl}
\nu(\{ M;\ M\models\psi(\iota)\})&=&
\nu_1(\{{\cal U};\ \{ n;\ \N\models\psi(n)\}\in \cal U\})\\
&=& dens_{u_0}\{ n;\ \N\models\psi(n)\}\\
&=& \frac 12,
\end{array}
\]
by 6. Thus the theorem is proved.
\end{proof}

\section{Remarks}\label{remarks}

1. We will generalize the concept of pseudorandom sequences and sets to
cover sequences in which $1$ occurs with frequency $p\neq\frac 12$. 

\bdf
Let $p$ be a real number, $0<p<1$. We will say that a sequence
$s:\N\to\{\pm 1\}$ is $p$-biased pseudorandom, if $s$ is computable in
nondeterministic polynomial time and 
\[
\lim_{n\to\infty}\sum_{i=0}^{n-1} f(i)((s(i)-1)/2 +p) =0,
\]
for every polynomial time computable function $f:\N\to\{\pm 1\}$. 

A set $X\sub\N$ will be called $p$-biased pseudorandom, if it is the set of
arguments for which a $p$-biased pseudorandom sequence is $1$.
\edf

Several propositions proved above generalize to $p$-biased
pseudorandom sequences and sets. In particular, we would like to draw
reader's attention to Corollary~\ref{cor1} and Theorem~\ref{t-main}. 
The condition 3. of Theorem~\ref{t-main} holds true for all real
numbers $p$, $0<p<1$ and all $p$-biased pseudorandom sequences
simultaneously. 

\bigskip
\noindent 2. The main weakness of Theorem~\ref{t-main} is that condition 3. is
stated only for one element, only for $\iota$. One can show that in
the constructed system $\cal M$, condition 3. holds for several other
elements of $M_0$; in particular, it is true for all elements of the
form $a\iota+b$ for $a\in\N$ and $b\in\Z$. This can easily be proved
using the following lemma.
\bl
Let $s$ be a pseudorandom sequence. Let $g\in\bf FP$ be increasing and
invertible in polynomial time on an infinite interval
$[n_0,\infty)$. Furthermore, suppose that $Rng(g)$, the range of $g$, has
positive density. Then the sequence $s'$ defined by $s'(x)=s(g(x))$ is
also pseudorandom.
\el
\bprf
Let $s$ and $g$ be given and suppose $s'$ is not pseudorandom. Let $f$
be a function that witnesses that $s'$ is not pseudorandom. We
define a function that witnesses that $s$ is not pseudorandom.
\[
f'(n)=\left\{
\begin{array}{ll}
f(g^{-1}(n)) &\mbox{ if }n\geq n_0\mbox{ and }n\in Rng(g),\\
0&\mbox{otherwise}.
\end{array}\right.
\]
\eprf
We certainly cannot expect condition 3. to hold for all numbers of
$M_0$. For small numbers $n$, $s(n)$ is defined in $M_0$, because $s$
is computable in exponential time. If $\alpha\in M_0$ is larger than
all numbers $c\iota$, for $c\in\N$, then $\alpha=[g]_{{\cal U}_0}$ for
some $g$ that grows more than linearly. The range of such a $g$ has
density 0, hence we cannot deduce anything about it. For example,
$\lambda(n^2)=1$ for all $n$, whence $M_0\models \lambda(\iota^2)=1$. 

\bigskip
\noindent 3. We also cannot expect stronger properties of random
sequences to hold in the system of models of Theorem~\ref{t-main}
unless we assume more about the sequences. For example, $s(\iota)$ and
$s(\iota+1)$ do not have to be independent in $\cal M$, because we do
not assume any kind of independence for pairs $s(n)$ and $s(n+1)$, 
$n\in\N$. A more specific example (assuming that $\lambda$ is pseudorandom)
is the fact that $\lambda(2n)=-\lambda(n)$.

\section{Philosophical speculations}

Some cosmologists believe that when the universe emerged from a
singularity some physical properties of it were decided randomly. 
Others even believe that there is a
{\em multiverse} consisting of many different universes, one of which
is our universe. In contrast to this, philosophers have never doubted
that the basic mathematical structures, namely the natural and real
numbers, are unique. These structures are unique in the sense that 
they must be same in all conceivable physical worlds.

There are good reasons to believe that the natural numbers are
absolute in the sense that there are no possible alternatives to them.
The only structures that satisfy the basic arithmetical laws and are
different from the natural numbers are nonstandard models.  Though it
has been proposed that the actual natural numbers have the structure
of a nonstandard model, e.g., in \cite{vopenka}, most philosophers do
not accept such a possibility. The problem is that a nonstandard model
contains the standard model as an initial part, and so we should
identify the natural numbers with this initial part. Thus viable
alternatives should use the same numbers with different arithmetical
operations. Since we can prove that the operations of addition and
multiplication are uniquely determined by the basic axioms (the
axioms of Robinson Arithmetic), it is inconsistent to assume that
on the set of standard numbers different kinds of addition and
multiplication are possible.

However, what is inconsistent in our world may be consistent in a
different one, and vice versa. Consider a pair of (necessarily
nonstandard) models of Peano Arithmetic $M$ and $N$ that have the same
elements, the same addition, but different multiplication. (Such pairs
can be easily constructed using recursively saturated models.)  In a
world in which $M$ is the standard natural numbers, it is inconsistent
to assume that anything like $N$ exists. Yet, it does.

We may secretly ponder over such scenarios, but there is a strong
reason not to talk about the possibility of different arithmetics
openly. If a concept is inconsistent, we cannot talk about it and
there cannot be any theory around it. 
Therefore, any conjecture of this kind would be neither provable nor
disprovable and, as such, should be discarded as meaningless. 

However, some phenomena can be studied even if they are not directly
observable---because they have side effects. The presence of these
effects is a proof of the phenomenon.  A side effect of the origin
of our integers in a random process could be the randomness present in
the structure of integers. It is not genuine randomness, because the
integers are a single structure. What we rather observe are some
properties that are satisfied by truly random objects. Therefore we
call it {\em pseudorandomness}. Number theorists are familiar with
this; they use assumptions about random behavior in heuristic
arguments when they are not able to prove theorems rigorously and some
conjectures are also justified in this way (including the Riemann
Hypothesis). 


\subsection*{Acknowledgment}
We would like to thank Emil Je\v{r}\'abek and Jan Kraj\'{i}\v{c}ek for
their useful comments and suggestions.


\begin{thebibliography}{99}

\bibitem{cook} S.A. Cook: Feasibly constructive proofs and the
  propositional calculus. In: Proc. seventh annual ACM symposium on
  Theory of computing, ACM New York, 83--97, (1975)

\bibitem{cook-nguyen} S.A. Cook, P. Nguyen: Logical Foundations of
  Proof Complexity. ASL series Perspectives in Logic, Cambridge
  Univ. Press, (2010)

\bibitem{goldreich} O. Goldreich: Foundation of Cryptography: Basic
  Tools. Cambridge Univ. Press, (2001)

\bibitem{green} B. Green: On (not) computing the M\"obius functions
  using bounded depth circuits. Combinatorics, Probability and
  Computing, to appear

\bibitem{kochen-kripke} S. Kochen, S. Kripke: Non-standard models of
  Peano arithmetic. In: Logic and arithmetic, int. Symp., Zürich 1980.
  Enseign. Math., II. S\'er. 28, 211--231, (1982)

\bibitem{krajicek} J. Kraj\'{\i}\v{c}ek: Bounded arithmetic,
  propositional logic, and complexity theory. Encyclopedia of
  Mathematics and Its Applications, Vol.60, Cambridge University
  Press, Cambridge - New York - Melbourne, (1995)

\bibitem{krajicek11} J. Kraj\'{\i}\v{c}ek: Forcing with random
  variables and proof complexity. London Mathematical Society Lecture
  Note Series, No.382, Cambridge University Press, (2011)

\bibitem{schnorr} C.P. Schnorr: A unified approach to the
  definition of a random sequence. Mathematical Systems Theory 5(3),
  246--258, (1971)

\bibitem{vopenka} P. Vop\v{e}nka: Mathematics in the Alternative Set
  Theory. Teubner, Leipzig (1979)

\end{thebibliography}
\end{document}